\RequirePackage{ifpdf}
\ifpdf 
\documentclass[pdftex]{sigma}
\else
\documentclass{sigma}
\fi

\begin{document}

\allowdisplaybreaks

\renewcommand{\PaperNumber}{022}

\renewcommand{\thefootnote}{$\star$}

\FirstPageHeading

\ShortArticleName{Laurent Polynomials and Superintegrable Maps}

\ArticleName{Laurent Polynomials and Superintegrable Maps\footnote{This paper is a contribution
to the Vadim Kuznetsov Memorial Issue ``Integrable Systems and
Related Topics''. The full collection is available at
\href{http://www.emis.de/journals/SIGMA/kuznetsov.html}{http://www.emis.de/journals/SIGMA/kuznetsov.html}}}

\Author{Andrew N.W. HONE}
\AuthorNameForHeading{A.N.W. Hone}

\Address{Institute of Mathematics, Statistics {\rm \&} Actuarial Science,\\
University of Kent, Canterbury CT2 7NF, UK}
\Email{\href{mailto:anwh@kent.ac.uk}{anwh@kent.ac.uk}}

\ArticleDates{Received October 26, 2006; Published online February 07, 2007}

\Abstract{This article is dedicated to the memory of Vadim Kuznetsov,
and begins
with some of the author's recollections of him.
Thereafter, a brief review of Somos sequences is provided,
with particular focus
being made on the integrable structure of Somos-4 recurrences, and on
the Laurent property. Subsequently a family of fourth-order recurrences
that share the Laurent property are considered, which are equivalent to
Poisson maps in four dimensions. Two of these maps turn out
to be superintegrable, and their iteration furnishes
inf\/initely many solutions of some associated quartic Diophantine equations.}

\Keywords{Laurent property; integrable maps; Somos sequences}

\Classification{11B37; 33E05; 37J35}

\renewcommand{\theequation}{\arabic{section}.\arabic{equation}}
\newcommand{\beq}{\begin{equation}}
\newcommand{\eeq}{\end{equation}}
\newcommand{\bea}{\begin{eqnarray}}
\newcommand{\eea}{\end{eqnarray}}
\newcommand{\bear}{\begin{array}}
\newcommand{\eear}{\end{array}}

\newcommand\mau{{\mathcal{u}}}
\newcommand\mav{{\mathcal{v}}}
\newcommand\ups{{\upsilon}}

\newcommand\LM{{\bf L}}
\newcommand\mM{{\bf M}}
\newcommand\mN{{\bf N}}
\newcommand\al{{\alpha}}
\newcommand\la{{\lambda}}
\newcommand\om{{\omega}}
\newcommand\uz{{\underline{z}}}
\newcommand\uk{{\underline{\kappa}}}
\newcommand\uu{{\underline{u}}}
\newcommand\uv{{\underline{v}}}

\newcommand\bu{{\bf {u}}}
\newcommand\bv{{\bf {v}}}
\newcommand\br{{\bf {r}}}
\newcommand\bx{{\bf {x}}}
\newcommand\by{{\bf {y}}}
\newcommand\bz{{\bf {z}}}

\newcommand\PS{{\mathcal P}}
\newcommand\larra{{\longrightarrow}}
\newcommand\TLM{{\bf \tilde{L}}}
\newcommand\TN{{\bf \tilde{N}}}
\newcommand\tq{{\tilde{q}}}
\newcommand\tp{{\tilde{p}}}

\renewcommand{\theequation}{\arabic{section}.\arabic{equation}}
\def\underset#1#2{\mathrel{\mathop{#2}\limits_{#1}}}
\newcommand{\haf}{{\hat{f}}}
\newcommand{\hax}{{\hat{x}}}
\newcommand{\hay}{{\hat{y}}}
\newcommand{\haz}{{\hat{z}}}
\newcommand{\hal}{{\hat{L}}}
\newcommand{\hak}{{\hat{K}}}

\newcommand\ka{{\kappa}}
\newcommand\be{{\beta}}
\newcommand\de{{\delta}}
\newcommand\gam{{\gamma}}
\newcommand\si{{\sigma}}
\newcommand\ze{{\zeta}}
\newcommand\epsi{{\epsilon}}

\newcommand\lax{{\bf L}}
\newcommand\mma{{\bf M}}
\newcommand\ctop{{\mathcal{T}}}
\newcommand\hop{{\mathcal{H}}}
\newcommand\ep{{\epsilon}}
\newcommand\T{{\tau}}
\newcommand\dd{\mathrm{d}}

\newcommand{\K}{{\mathbb K}}  
\newcommand{\F}{{\mathbb F}}  
\newcommand{\Pro}{{\mathbb P}}  
\newcommand{\R}{{\mathbb R}}
\newcommand{\N}{{\mathbb N}}

\newcommand{\Q}{{\mathbb Q}}
\newcommand{\Z}{{\mathbb Z}}
\newcommand{\C}{{\mathbb C}}
\newcommand{\sh}{\mathcal{S}}

\newcommand{\vN}{{\mathcal{N}}}
\newcommand{\J}{{\mathcal{J}}}
\newcommand{\vL}{{\mathcal{L}}}
\newcommand{\vM}{{\mathcal{M}}}

\newcommand\Lop{{\mathrm{L}}}

\newcommand{\rd}{{\mathrm{d}}}
\newcommand{\sgn}{{\mathrm{sgn}}}
\def \ring {{\cal R}}

\newcommand{\B}{{\mathcal{B}}}
\newcommand{\A}{{\mathcal{A}}}
\newcommand{\Cs}{{\mathcal{C}}}
\newcommand{\rs}{{\mathcal{S}}}
\newcommand{\rA}{{\mathcal{A}}}
\newcommand{\rB}{{\mathcal{B}}}
\newcommand{\tilta}{{\tilde{\tau} }}
\newcommand{\tiJ}{{\tilde{\mathcal{J}}}}
\newcommand{\rQ}{{\mathcal{Q}}}

\renewcommand{\thefootnote}{\arabic{footnote}}
\setcounter{footnote}{0}

\section{Introduction}

It is with considerable sadness that I begin to write this piece in memory of Vadim
Kuznetsov, whose death came as a great shock to me.
However, I do not wish to remain in melancholic
mode,
but rather I would like to recall some of my fondest and happiest
memories of him.

While I was a PhD student in Edinburgh, I used to travel
to Leeds every so often to attend the LMS workshops on integrable systems,
and I'm sure I must have f\/irst met Vadim at one of these meetings.
To begin with I remember his charming smile, as well as his relaxed way of
asking penetrating mathematical questions. I also recall the great
enthusiasm and energy with which he would give a seminar, and his clarity of
presentation.

Shortly after graduating from Edinburgh, in September 1997
I went to Rome to take up my f\/irst postdoctoral
position, working with Orlando Ragnisco in the Physics Department of Roma Tre.
It was during this period that I had the privilege of getting to know Vadim
a lot better. A~few months after my arrival, he came to Rome to visit
Orlando for a month, and the three of us ended up working together on a
project that was suggested by Vadim,
concerning B\"acklund transformations (BTs) for f\/inite-dimensional integrable
Hamiltonian systems. This turned out to be very fruitful,
resulting in three joint publications \cite{hkr1,hkr2,hkr3}.

Vadim's presence in Rome was immensely stimulating for me, because
he succeeded in posing just the right question, at the
precise moment when I had the necessary tools available to answer~it.
The specif\/ic problem that he f\/irst presented to me and Orlando was the construction
of BTs for certain integrable classical mechanical systems corresponding to reduced
Gaudin magnets.
A~particular concrete example of such a system was the case (ii) H\'enon--Heiles system,
an integrable system with two degrees of freedom. As it happened, in my PhD thesis I had
already constructed an analogous BT for the non-autonomous case of this
system, as well as deriving the explicit formula for the generating function
of the canonical (contact) transformation in that case \cite{nahh}. During
my viva voce examination a few months earlier, Allan Fordy
had actually asked me  whether the same sort of derivation could also be applied to
the autonomous case, to produce a Poisson correspondence
in the spirit of \cite{poissonmaps}, and I could see no obvious
obstruction.
Thus it was that,
when Vadim arrived in Rome,
his vivid explanation of BTs, as well as his insistence that we should
start constructing new ones, was all that I needed to work out the BT for the
H\'enon--Heiles system \cite{hkr1}, and this soon revealed a similar algebraic structure
underlying many other examples \cite{hkr2}.

After he left Rome, I saw Vadim again in June 1998 at the conference
\textit{Integrable Systems: Solutions and Transformations} in Guardamar, Spain, where
he came with his wife, Olga, and his son, Simon. We sat down together in the sunshine and completed some
of the work on the second paper \cite{hkr2} while we were there. Subsequently,
I saw Vadim sporadically at various meetings in Leeds and elsewhere, and we always found
the time for a friendly chat about our lives and work. I particularly remember a
very brief and enjoyable
(but f\/iery)
dispute that we had in Cambridge in 2001, while
sitting together during an interlude between lectures in the Newton Institute. It boiled
down to a minor dif\/ference in our points of view, which we respectively argued for without
compromise, so that (having each seen the other's perspective) there was no love lost between us.

The rest of this article is concerned with a family of discrete dynamical systems (Poisson maps) in four
dimensions, the f\/irst few of which are integrable, while the others are not. Before
going into details, I should like to explain why I have chosen this topic. The work I did in
my PhD was primarily concerned with  integrable systems in the continuous setting (ordinary and
partial dif\/ferential equations), and it was not until Vadim's visit to Rome that I began
to get actively interested in discrete systems. Ever since then, I have found the subject of
discrete dynamics increasingly fascinating, and I shall always have Vadim to thank for
inspiring me to look in this direction. Another interesting and unexpected property of the
Poisson maps considered below is that their iterates are Laurent polynomials in the initial data;
this is
an instance of the Laurent phenomenon \cite{fz}. Vadim was an expert on special functions, and orthogonal
polynomials in particular (for one of his many contributions in this area, see
\cite{qopsjack}, for instance). However, most of the sequences of (Laurent) polynomials
treated below satisfy nonlinear equations instead of linear ones.

The theory of discrete integrable maps has seen a great deal of activity in the past twenty
years. The situation was much clarif\/ied by Veselov \cite{ves1,ves2,ves3} who introduced
integrable Lagrange correspondences -- a natural discrete-time analogue of Liouville integrable
continuous f\/lows -- which induce (generically multi-valued) shifts on the associated Liouville
tori (see also~\cite{rag}).
Given a continuous integrable system, it is natural to seek a discretization
of it that retains both the integrability and as many other properties as possible (e.g.\ Poisson
structure, Lax pair, etc.). However, in general such a time-discretization will be implicit,
and it will not preserve the same integrals as the original continuous system (see
\cite{surisbook} for the state of the art in integrable discretizations). Building on results
obtained for the Toda lattice by
Pasquier and Gaudin \cite{pg}, Kuznetsov and Sklyanin identif\/ied a special
class of time-discretizations for
integrable Hamiltonian systems which they referred to as BTs \cite{kuskly}, by analogy with
B\"acklund transformations for evolutionary PDEs.

In the  setting of
f\/inite-dimensional systems with a Lax pair,
BTs were identif\/ied as explicit Poisson maps which preserve the same set of integrals as
the continuous f\/low that they discretize, and depend on a B\"acklund parameter $\la$ which
satisf\/ies a certain `spectrality' property with respect to a conjugate variable $\mu$ (where
$(\la , \mu )$ are the coordinates of a point on the spectral curve associated with the Lax pair).
The viewpoint that I emphasized in \cite{hkr1,hkr2,pinney} was that the systems being considered
were reduced/stationary f\/lows of the KdV hierarchy, whose BTs could be obtained
by reduction from the Darboux--B\"acklund transformation for KdV (this is in the same vein
as the dressing chain \cite{dressing} -- see also \cite{weiss}), while the BTs in \cite{hkr3} were
derived more directly. In Vadim's work with Pol Vanhaecke \cite{vadimpol},
all of the previously known examples were unif\/ied via an algebro-geometric approach,
which explained the deeper meaning of BTs as discrete shifts on the (generalized) Jacobian of
the associated spectral curve,
thus identifying them as the discrete-time counterparts
of algebraically completely integrable systems,
as described in \cite{vanhaecke}, for instance.
While there has been subsequent work by Vadim and others
on BTs in classical mechanics \cite{chm,fed2,kupera},
a lot of the original motivation for studying them  came
from {\it quantum} integrable  systems (Baxter's $Q$-operator). This
idea has proved extremely ef\/fective (see e.g.~\cite{kuskly2,qopsjack}), and will no doubt
continue to bear fruit for a long time to come.

The last time I saw Vadim was in Leeds in April 2005, when he invited me to give one
in the series of Quantum Computational seminars that he organized
there\footnote{\href{http://www.maths.leeds.ac.uk/~vadim/QCS.htm}{http://www.maths.leeds.ac.uk/\~{}vadim/QCS.htm}}.
At that time I~spoke about Somos sequences, which are reviewed
in the next section. In the evening
after the seminar I went out for a very enjoyable dinner with Vadim and Olga, together with Oleg Chalykh
and Sara Lombardo.  I made an appointment to see Vadim in his of\/f\/ice early the next morning, so that before
my return home we
had a good discussion about his recent work on the integrable dynamics of spin chains that arise in
models of
Fermi--Bose condensates \cite{yka} and BCS superconductors \cite{yake}, and he described an unsolved
problem concerning special solutions.   This is how I remember him now: full of energy and always
seeking to answer new questions.

\section{Somos sequences and the Laurent property}

The properties of integer sequences generated by linear recurrences have been the subject
of a~great deal of study in number theory, and nowadays they f\/ind applications
in computer science and cryptography \cite{recs}. However, the theory of nonlinear recurrence
sequences is still in its infancy. Clearly, a $k$th-order
nonlinear recurrence relation of the form
\beq \label{rec1}
x_{n+k}=F(x_n,x_{n+1},\ldots ,x_{n+k-1})
\eeq
is just a particular sort of discrete dynamical system, so such recurrences can be considered as
generating a special type of nonlinear dynamics. If we want (\ref{rec1}) to generate sequences of
integers, then choosing $F$ to be a polynomial with integer coef\/f\/icients will certainly do the trick,
but
in general the corresponding
map in $\R^k$ (or $\C^k$) will not have a unique inverse.
Moreover, in that case such sequences
generically exhibit
double exponential growth i.e.\ $\log |x_n|$ grows exponentially with $n$. A simple example in this class
is the
quadratic map def\/ined by the recurrence
\beq \label{logi}
x_{n+1}=x_n^2+c
\eeq
with a parameter $c$,
which is a prototypical model of chaos. However, note that
the special cases $c=0,-2$
are exactly solvable \cite{chaos}, and in these cases
one can also argue
that (\ref{logi}) is integrable in the sense of admitting a commuting map
(see \cite{ves1} and references).
The theory of linear recurrence sequences relies heavily on the fact that they are explicitly solvable.
Thus it is natural to look for nonlinear recurrences that share this
property, or that are integrable in a~broader sense.

In the case that the map corresponding to (\ref{rec1})
is invertible,
one can  also
allow $F$ to be a~rational function, thereby
considering birational maps, but then it is no longer
clear that integer sequences should result. However, it turns out that among 
those rational recurrences of the particular form
\begin{gather} \label{rec2}
x_{n+k}\, x_n=f(x_{n+1},\ldots ,x_{n+k-1}),
\end{gather}
there is a very large class of recurrences that generate integer sequences
from suitable initial data.
One of the f\/irst known examples of this type is the Somos-4 recurrence
\beq \label{s4}
x_{n+4}\,x_n = \al \, x_{n+3}\, x_{n+1} + \be \, (x_{n+2})^2,
\eeq
which was found by Michael Somos when he was investigating the combinatorics
of elliptic theta functions. Somos observed numerically that by taking
the coef\/f\/icients $\al =\be =1$ and initial data $x_0=x_1=x_2=x_3=1$,
the fourth-order recurrence (\ref{s4})
yields a sequence of integers \cite{sloane}, that is
\beq \label{s4seq}
1,1,1,1,2,3,7,23,59,314,1529,8209, 83313, \ldots.
\eeq
Similarly he noticed that for the Somos-$k$ recurrences
\beq \label{sk}
x_{n+k}\, x_n
=\sum_{j=1}^{[k/2]} \al_j\, x_{n+k-j}\,x_{n+j}
\eeq
with all coef\/f\/icients $\al_j=1$, if all $k$ initial values are 1 then
an integer sequence results for $k=5,6,7$, but denominators appear for $k=8$.

Various direct proofs that the terms of the sequence (\ref{s4seq})
are all integers were found at the beginning of the  1990s,
when various other examples were found \cite{gale, rob}, but
a deeper understanding came from the realization that the recurrence (\ref{s4}) has the
{\it Laurent property}: its iterates are all Laurent polynomials in the initial data
(and in $\al$, $\be$)
with integer coef\/f\/icients. To be more precise, the iterates
of (\ref{s4}) satisfy $x_n\in \Z [x_0^{\pm 1}, x_1^{\pm 1}, x_2^{\pm 1}, x_3^{\pm 1}, \al ,\be ]$
for all $n$, from which the integrality of the particular sequence (\ref{s4seq}) follows
immediately. A little earlier, when
Mills, Robbins and Rumsey made their study of the Dodgson condensation method for
computing determinants \cite{mrr} (which produced the famous alternating sign matrix
conjecture \cite{bressoud}), they considered the recurrence
\beq
D_{\ell ,m,n+1}D_{\ell ,m,n-1}=\al \, D_{\ell +1,m,n}D_{\ell -1,m,n} +\be \,
D_{\ell ,m+1,n}D_{\ell ,m-1,n} ,
\label{hirota}
\eeq
for $\al =1$ and
observed that it produced Laurent polynomials in the initial data.
The equa\-tion~(\ref{hirota}) thus
became known within the algebraic combinatorics community, where it is referred to as the
octahedron recurrence \cite{propp},
while in the theory of integrable systems it is known as a particular form of the
discrete Hirota equation \cite{zabrodin} (the bilinear equation for
the tau-function of discrete KP). The
Somos-4 recurrence (\ref{s4}) is an ordinary dif\/ference reduction of the
partial dif\/ference equation (\ref{hirota}): it has been noted by Propp
that if $x_n$ satisf\/ies (\ref{s4}) then $D_{\ell ,m,n}=x_{2n+m}$ satisf\/ies
the discrete Hirota equation (see also \cite{speyer} for another reduction).

Many more examples of this Laurent property have begun to emerge
quite recently as an of\/fshoot of the theory of cluster
algebras due to Fomin and Zelevinsky (see \cite{fz4} and references).
The exchange relations in a cluster algebra of rank $k$ are typif\/ied by a recurrence of the form
\beq\label{cluster}
x_{n+k}\,x_n=c_1\, M_1 (x_{n+1},\ldots ,x_{n+k-1}) + c_2\, M_2 (x_{n+1},\ldots ,x_{n+k-1})
\eeq
for suitable monomials $M_j$ and coef\/f\/icients $c_j$, which is a special case
of (\ref{rec2}).
In \cite{fz}, the general machinery of cluster algebras was shown to be very ef\/fective
in proving the Laurent property for a wide variety of recurrences, mostly
(but not all) of the
form (\ref{cluster}).
In particular, Fomin and Zelevinsky there gave the f\/irst proof of the Laurent property for
the octahedron (discrete Hirota) recurrence
(\ref{hirota}). Subsequently, Speyer has developed a combinatorial model to prove
more detailed properties of the Laurent polynomials generated by this recurrence --
in particular, that all the coef\/f\/icients are 1 \cite{speyer}.

So far we have discussed the integrality of the sequence (\ref{s4seq}), but not
the integrability of the Somos-4 recurrence. Taking
$(x_0,x_1,x_2,x_3)$ as coordinates, the map $\C^4 \to \C^4$
corresponding  to~(\ref{s4}) preserves the degenerate Poisson bracket def\/ined by
\beq\label{s4brkt}
\{ \, x_m , x_n \, \}_2 =(n-m)\, x_mx_n,
\eeq
which has Casimirs
\beq \label{un}
u_n=\frac{x_{n-1}x_{n+1}}{(x_n)^2 }, \qquad n=1,2.
\eeq
This bracket is of the `log-canonical' type that has previously been found
in the context of cluster algebras \cite{gekhtman}; it is natural to
consider it as a Poisson bracket on the f\/ield of rational functions
$\C (x_0,x_1,x_2,x_3)$. (The reason for the subscript $2$ on the bracket will
become apparent in the next section.)
The set of solutions of (\ref{s4}) is invariant under the two-parameter
Abelian group of gauge transformations generated by
\beq \label{gauge}
x_n \mapsto A\, x_n, \qquad x_n \mapsto  B^n \, x_n ,\qquad A,B\in \C^*.
\eeq
The Hamiltonian vector f\/ields corresponding to these transformations
are respectively generated by the rational monomials
\beq \label{gaugeflows}
K_1 = \frac{x_1}{x_2}, \qquad K_2=\frac{(x_1)^2}{x_2},
\eeq
which satisfy
\[
\{\, K_1,x_n \,\}_2 =K_1 \, x_n, \qquad
\{\, K_2,x_n \,\}_2 =n \, K_2\, x_n, \qquad
\{ \, K_1,K_2\, \}_2=K_1K_2.
\]

In fact, the most interesting part of the dynamics generated by  (\ref{s4})
takes place in the plane spanned by the Casimirs $u_1$, $u_2$ for the bracket $\{\cdot ,\cdot\}_2$. If we take
the def\/inition  (\ref{un}) to hold for all $n$, then the quantities $u_n$ are
clearly invariant under the gauge transformations (\ref{gauge}), and
satisfy the second-order recurrence
\beq
\label{ueq}
u_{n+2}=\frac{\al u_{n+1} +\be }{u_n(u_{n+1})^2} .
\eeq
(So the fourth-order equation (\ref{s4}) is the
Hirota bilinearization of (\ref{ueq}), which is
a second-order ordinary dif\/ference equation.)
By taking
$(u_1,u_2)$ as coordinates in  $\C^2$, this corresponds to the rational map
of the plane given by
\[
\left(\bear{c} u_1 \\ u_2 \eear \right) \mapsto
\left(\bear{c} u_2 \\ (\al u_{2} +\be )/(u_1(u_{2})^2) \eear \right)
\]
which
preserves the Poisson bracket
\beq \label{newpb}
\{ \, u_n,u_{n+1}\, \} =u_n u_{n+1},
\eeq
or equivalently the symplectic form
\beq \label{om}
\om_n = (u_nu_{n+1})^{-1} \, \dd u_n \wedge \dd u_{n+1}
\eeq
such that $\om_{n+1}=\om_n$. Furthermore, this
has the conserved quantity
\beq \label{J}
J=u_n\, u_{n+1}
+\al \,\left(\frac{1}{u_n}+\frac{1}{u_{n+1}} \right)
+\frac{\be}{u_n \, u_{n+1}},
\eeq
which def\/ines a quartic curve
\beq \label{quartic}
\ze^2 \eta^2 -J\, \ze \eta +\al (\ze + \eta ) +\be =0
\eeq
of genus one. Hence we see that (\ref{ueq}) produces a Liouville integrable
system with one degree of freedom, and the curve (\ref{quartic})
itself def\/ines the two-valued correspondence $u_n \mapsto u_{n\pm 1}$, which is
a particular case of the Euler--Chasles  correspondence (see \cite{ves1,ves2,ves3}).

Upon uniformizing the elliptic quartic we f\/ind that the explicit solution
to (\ref{ueq}) is given by
\beq\label{usol}
u_n = \wp (z) - \wp (z_0 + nz),
\eeq
in terms of the Weierstrass $\wp$ function for the elliptic curve
\beq\label{weier}
E: \quad Y^2=4X^3-g_2X-g_3,
\eeq
with
$
g_2=12\la^2-2J$, 
$g_3=4\la^3-g_2\la-\al$,
$\la = (J^2/4-\be )/(3\al )$,
%
and  $z_0$, $z\in \C /\Lambda =\mathrm{Jac}(E)$ are given by elliptic integrals
obtained from inversion of the relations
$
\wp (z) =\la $, 
$\wp (z_0 )= \la - u_0 $.  The coef\/f\/icients
$\al$, $\be $ and also $J$
are given as elliptic functions of $z$ by 
$\al =\wp'(z)^2$, $\be = \wp'(z)^2\, (\wp (2z)-\wp (z))$, $J=\wp''(z)$.

From this it follows \cite{honeblms}
that the solution to the initial value problem
for the Somos-4 recurrence (\ref{s4}) can be written in terms
of the Weierstrass sigma function as
\beq \label{s4sol}
x_n =A B^n \,
\frac
{\si (z_0 + n z)}
{\si (z) ^{n^2}}
\eeq
for suitable $A$, $B$. There is an analogous formula for the
general solution of the Somos-5 recurrence \cite{hones5}, which has an
additional dependence on the parity of $n$.


\begin{figure}[ht!]
\centerline{\scalebox{0.4}{\includegraphics[angle=270]{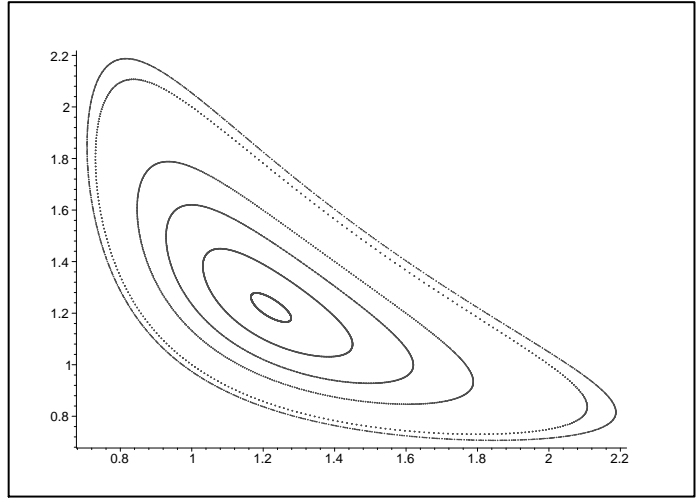}}}
\caption{A family of orbits for the nonlinear
recurrence (\ref{ueq}) associated with Somos-4.}
\label{s4orbits}
\end{figure}

The map def\/ined by (\ref{ueq}) is a very simple example of the QRT family \cite{qrt}.
It has a $2\times 2$ discrete Lax pair given by
\beq \label{lax}
\LM_n \, \Psi_n = \nu\,\Psi_n, \qquad  \Psi_{n+1} = \mM_n \, \Psi_n,
\eeq
where
\[
\LM_n (\ups ) = \left( \bear{ccc} u_n u_{n+1} && -\al \\

\ups -u_n-u_{n+1} && \al (1/u_n +1/u_{n+1} ) +\be /(u_n u_{n+1}) \eear \right),
\]
\[
\mM_n (\ups )  = \left(\bear{ccc} 0 && -\al \\

\ups -u_n-u_{n+1} && \al /u_n  \eear \right).
\]
The equation (\ref{ueq}) arises as the compatibility condition
$\LM_{n+1}\mM_n =\mM_n\LM_n$ for the system~(\ref{lax}), and the
associated spectral curve
is given by
\beq\label{spec}
\mathrm{det}\, (\LM_n (\ups ) -\nu \, \boldsymbol{1} ) =
\nu^2-J\nu +\al\ups +\be=0 ,
\eeq
where $J$ is the conserved quantity given by (\ref{J}).
From the formulae
$
\nu = \ze \eta$,
$\ups =\ze  +\eta
$
the elliptic quartic curve (\ref{quartic}) is seen to be a ramif\/ied double
cover of this rational (genus zero) spectral curve.

It is clear from the above considerations that the dynamics of (\ref{ueq})
corresponds  to a sequence of points $P_0 + nP$ on the elliptic curve
$E$ given by (\ref{weier}), or to the equivalent discrete linear f\/low
$z_0 + n z$ on its Jacobian, and in that sense (as was noted in \cite{honeblms})
it is the same as the underlying dynamics of the BT for the one-particle
Garnier system  constructed in \cite{hkr2}, or that of the BT for the $g=1$ odd Mumford
system as in
\cite{vadimpol}. However, while one can make changes of variables
between (\ref{ueq}) and each of the latter two BTs, they are not canonical
transformations, because the Poisson bracket (\ref{newpb})
is incompatible with the Poisson  structures of either of these BTs.
Nevertheless, just as for the BTs, the recurrence (\ref{ueq}) is a discretization of
a continuous time integrable system with the same Poisson structure and
conserved quantities (in this case, only one of them), namely the
f\/low in the plane with Hamiltonian $J$ def\/ined by (\ref{J}) with $n=1$,~i.e.
\begin{gather}
\frac{\dd u_1 }{\dd t}  = \{ \, J,u_1\,\}  =  (\al u_1+\be )/u_2-(u_1)^2u_2,
\nonumber\\
\frac{\dd u_2 }{\dd t}  = \{ \, J,u_2\,\}  =
-(\al u_2+\be )/u_1+(u_2)^2u_1.\label{jsys}
\end{gather}
From the same uniformization of the quartic (\ref{quartic}) as before,
the solution of the  system (\ref{jsys}) can be written down as
\[
u_1(t)=\wp (z) -\wp (\tilde{z}_0+\wp'(z) t),
\qquad
u_2(t)=\wp (z) -\wp (\tilde{z}_0+z+\wp'(z) t),
\]
and upon comparison with (\ref{usol}) it can be seen directly
how the discrete f\/low interpolates the continuous one
(cf. Fig.~\ref{s4orbits}).

The construction of a sequence of points $P_0 + nP$ on elliptic curve $E$
from
a Somos-4 or \mbox{Somos-5} sequence was previously understood in
unpublished work of several number theorists\footnote{For further references see
\href{http://www.math.wisc.edu/~propp/somos.html}{http://www.math.wisc.edu/\~{}propp/somos.html}.}~--
see the discussion of Zagier \cite{zagier},
and the results of Elkies quoted in \cite{heron}.
The  algebraic part of the construction is described in the thesis of
Swart \cite{swart} (who also mentions unpublished results of Nelson
Stephens), and van der Poorten has recently presented another
construction based on the continued fraction expansion of the square root of
a quartic \cite{vdp}. In fact,
Somos-4 sequences 
have an ancestor from the 1940s,
in Morgan Ward's  work on elliptic divisibility sequences (EDS),
which just correspond
to multiples of a point $n\,P\in E$  \cite{ward1, ward2}
i.e. this is the special case $P_0 = \infty$, so that
$z_0=0$, with the further requirement that
$A=B=1$ in (\ref{s4sol}).
The iterates of an EDS, which are generated by
(\ref{s4}) with coef\/f\/icients $\al =(x_2)^2$, $\be =-x_1x_3$ and integer initial
data  $x_1=1$, $x_2,x_3,x_4\in\Z$ with $x_2|x_4$, satisfy the divisibility property
$x_m|x_n$ whenever $m|n$, and correspond to values of the division
polynomials of the curve (for a description of these see Exercise 3.7 in \cite{silver1}).
In this  sense, an EDS generalizes properties of certain linear recurrence sequences.
For example, the
Fibonacci numbers are generated by the recurrence $F_{n+1}=F_n+F_{n-1}$ with initial values
$F_0=1$, $F_1=1$, and form a divisibility sequence. Moreover, the even index terms
$x_n=F_{2n}$  form a divisibility sequence (so $F_{2m}|F_{2n}$ whenever $m|n$)
and also satisfy the Somos-4 recurrence
\[
F_{2n+4}\,F_{2n-4} = 9F_{2n+2}\,F_{2n-2} -8(F_{2n})^2 ,
\]
which corresponds to a degenerate case of the curve (\ref{weier})
where the discriminant vanishes, so $g_2^3-27g_3^2=0$ and the
formula (\ref{s4sol}) degenerates to an expression in terms of the hyperbolic sine.

The arithmetical
properties of EDS and Somos sequences -- in particular the distribution
of primes therein -- are a subject of current
interest \cite{eew, ems, silverman}. Some of these properties 
are discussed in the book \cite{recs} (see section 1.1.20, for instance),
where it is suggested that such bilinear recurrences should be suitable
generalizations of linear ones, with many analogous features. Based on the
appearance of higher-order Somos recurrences in the work of Cantor
on the analogues of division polynomials for hyperelliptic curves \cite{cantor}
(see also \cite{matsupsi} for analytic formulae), it was conjectured in \cite{honeblms}
that every Somos-$k$ sequence should correspond to a discrete linear
f\/low on the Jacobian of such a curve (with an associated discrete integrable system),
and the plausibility of this conjecture was justif\/ied
by a na\"ive counting argument. However, on Propp's bilinear forum \cite{proppweb},
Elkies had already given a more detailed
argument to the contrary,  based on a proposed theta function formula for the terms of such
sequences, which indicated that while the general Somos-6 and Somos-7 sequence
could be described by such a formula in genus two, the general Somos-$k$ for $k\geq 8$
could not. Thus in this setting the absence of the Laurent property appears to coincide
with the absence of algebraic integrability.

Nevertheless, in \cite{beh} it was shown that a particular
family of solutions of Somos-8 recurrences
can be described in terms of the Kleinian sigma
function for a genus two curve (which is equivalent to an expression
in theta functions), and these solutions are related to the BT for the
H\'{e}non--Heiles system that was found in \cite{hkr1, hkr2}.
The author has also found that the Somos-6 and Somos-7 recurrences
correspond to a rational map in $\C^4$ with two independent conserved
quantities, and there is a similar expression for the solutions in terms
of genus two sigma functions. For instance,
letting $\si $ denote the genus two Kleinian
sigma function (see e.g.~\cite{bel} for the def\/inition),
associated with a curve given by the af\/f\/ine
equation $y^2 = 4x^5 +c_3 x^3+\dots + c_0$ with period lattice $\Lambda$,
the expression
\beq \label{s6sol}
x_n = A\, B^n \, \frac{\si ({\bv}_0 + n\, \bv )}{\si (\bv )^{n^2}},
\eeq
where $A,B \in \C^*$, ${\bv}_0 , \bv \in \C^2\bmod \Lambda$,
satisf\/ies a Somos-6 recurrence provided
that $\bv$ is constrained according to
\[
\left| \bear{ccc}
1 & 1 & 1 \\
\wp_{12}(\bv ) & \wp_{12}(2\bv ) &\wp_{12}(3\bv )  \\
\wp_{22}(\bv ) & \wp_{22}(2\bv ) &\wp_{22}(3\bv )
\eear
\right| =0,
\]
where $\wp_{jk} (\bv )= -\partial_j\partial_k \log \si (\bv )$
are the corresponding Kleinian $\wp$ functions. In the case of generic~$\bv$,
if this  constraint does not hold, then $x_n$ given
by (\ref{s6sol}) satisf\/ies a Somos-8 recurrence instead. The full
details will be presented elsewhere~\cite{ddiis}.

Before moving onto other examples in the next section, we should mention
one more feature of the Somos-4 recurrence, namely the fact that it generates
solutions of a quartic Diophantine equation in four variables. If we rewrite the formula
(\ref{J})
for the conserved quantity $J$ in terms of the original variables $x_n$,
we obtain the
equation
\begin{gather}
(x_{n-1})^2 (x_{n+2})^2+\al \, \Big(x_{n-1}(x_{n+1})^3+ (x_{n})^3x_{n+2}
\Big) + \be \,  (x_{n})^2 (x_{n+1})^2 \nonumber\\
\qquad{}
=  J\,x_{n-1}x_{n}x_{n+1}x_{n+2} .\label{s4dio}
\end{gather}
If we have coef\/f\/icients $\al , \be \in \Z$ (or in $\Q$),
and if the Somos-4 recurrence (\ref{s4})
with a set of integer initial data $(x_0,x_1,x_2,x_3)$
generates a non-periodic
sequence of iterates satisfying $x_n\in \Z$ for all $n$,
then there are inf\/initely many quadruples of integers
$(x_{n-1},x_n,x_{n+1},x_{n+2})$ that are solutions
of the quartic Diophantine equation (\ref{s4dio}). (Note that
in this case,
as long as all the integer initial data are non-zero, then the coef\/f\/icient $J$
which appears in (\ref{s4dio}) is uniquely determined, and $J\in\Q$.) This
can be seen as a particular instance of a general feature shared by
all  recurrences that both have the Laurent property and possess a rational
invariant: generically, the orbit of suitable initial data will  generate
inf\/initely many solutions of an associated Diophantine equation.

\medskip
\noindent {\bf Diophantine Laurentness Lemma.\/}
{\it Suppose that a $k$th-order rational recurrence of the form~\eqref{rec1}
has coefficients in $\Q [{\bf c}]$ (for
some set of parameters ${\bf c}$) and has the Laurent pro\-per\-ty,
i.e.\ $x_n\in \Z [x_0^{\pm 1},x_{1}^{\pm 1},\ldots ,x_{k-1}^{\pm 1},{\bf c}]$ for all $n$.
Suppose further that this recurrence also has a rational conserved quantity given by
\[
K=\frac{f_1(x_n,x_{n+1},\ldots ,x_{n+k-1},{\bf c})}
{f_2(x_n,x_{n+1},\ldots ,x_{n+k-1},{\bf c})}
\]
for $f_1,f_2\in \Z [x_n,x_{n+1},\ldots ,x_{n+k-1},{\bf c}]$.
If $f_2\neq 0$ for some fixed integer values of ${\bf c}$ and initial
data $x_j = 1$ or $-1$ for $j=0,\ldots, k-1$, then the
value of $K\in\Q$ is fixed, and the recurrence
generates infinitely many integer
solutions of the Diophantine equation
\[
f_1(x_n,x_{n+1},\ldots ,x_{n+k-1},{\bf c})=K\, f_2(x_n,x_{n+1},\ldots ,x_{n+k-1},{\bf c})
\]
as long as the
corresponding orbit is not periodic.
}

\medskip

The integer sequence (\ref{s4seq}) provides a concrete example of the above
result: setting $\al =\be =1$, the initial data $1,1,1,1$ yield the value
$J=4$ in (\ref{s4dio}), and for $n\geq 0$ any four adjacent
terms of this increasing sequence provide a distinct solution of the equation.
In fact, in \cite{swahon}  it is proved that the
iterates of the Somos-4 recurrence satisfy the stronger property that
$x_n\in \Z [x_0^{\pm 1}, x_1,x_2,x_3,\al ,\be ,(\al^2 +\be J) ]$ for $n\geq 0$,
which yields a broader set of suf\/f\/icient criteria for integer sequences
to be generated.
In the next section we will see
analogous results for some other recurrences.

\section{A fourth-order family}

\setcounter{equation}{0}


A 
generalization of (\ref{s4}), that retains the Laurent
property,
is
the family of fourth-order recurrences
\begin{gather} \label{7m123}
x_{n+2}\, x_{n-2} = x_{n+1}^{a}\, x_{n-1}^{b}  +  x_n^{c},
\end{gather}
where the exponents $a$, $b$, $c$ are positive integers. These generalized
Somos-4 recurrences were
f\/irst described in print by David Gale \cite{gale}, who noted that they
all generate integer sequences from the
initial values $x_0=x_1=x_2=x_3=1$. Among various examples
covered in \cite{fz}, Fomin and Zelevinksy subsequently
proved that each of these recurrences has the Laurent
property. However, in contrast to the integrable structure
of the original Somos-4 recurrence, most of these examples
do not seem to correspond to  completely integrable systems.

Below we shall not present an analysis
of the complete family (\ref{7m123}),
but rather we focus on the special sub-family
of recurrences
def\/ined by $a=b=1$, with
$c\in\N$. In this case, it will be convenient to introduce a
parameter $\be$ as
the coef\/f\/icient of the second term on the right hand side;
although this can always be removed by rescaling $x_n$, its inclusion
preserves the Laurent property (while inserting another coef\/f\/icient $\al$
in front of the $x_{n+3}x_{n+1}$ term does not, unless $c=2$). These recurrences
also satisfy the singularity conf\/inement test that was proposed in \cite{grp} as
an analogue of the Painlev\'e test for discrete equations:
if an apparent singularity is reached (in this case, corresponding
to the situation that one of the iterates vanishes), then it is always possible
to analytically continue through it.

\begin{proposition}
For each $c\in\N$ the recurrence
\beq \label{7betarec}
x_{n+4}\, x_{n} = x_{n+3}\, x_{n+1} +  \be \, x_{n+2}^c,
\eeq
which corresponds to the iteration of the rational map
\beq \label{7m123map}
\left( \begin{array}{c} x_0 \\x_1 \\ x_2 \\ x_3 \end{array} \right) \mapsto
\left( \begin{array}{c} x_1 \\x_2 \\ x_3 \\ (x_1x_3 + \be\, x_2^{c})/x_0
\end{array} \right) ,
\eeq
has the Laurent property
in the sense that
$x_n\in\Z [x_0^{\pm 1},x_{1}^{\pm 1},x_{2}^{\pm 1},x_{3}^{\pm 1},\be ]$ for
all $n\in\Z$, and also satisfies the singularity confinement test.
Furthermore,  \eqref{7m123map} is a Poisson map with respect to
the   log-canonical Poisson bracket $\{\cdot ,\cdot \}_c$ defined by
\begin{alignat}{4}
&\{\, x_0,x_1\,\}_c=x_0x_1 ,\qquad && \{\, x_0,x_2\,\}_c=cx_0x_2 ,\qquad && \{\, x_0,x_3\,\}_c=(c+1)x_0x_3,&\nonumber  \\
& \{\, x_1,x_2\,\}_c=x_1x_2, && \{\, x_1,x_3\,\}_c=cx_1x_3, &&\{\, x_2,x_3\,\}_c=x_2x_3 ,&\label{pbs}
\end{alignat}
which is nondegenerate for $c\neq 2$.
\end{proposition}

\begin{proof}
The recurrence (\ref{7betarec}) is of the cluster algebra type,
so the Laurent property can be proved
by the methods of \cite{fz}, where the details for the complete family
(\ref{7m123}) are presented. However, here it is convenient to
sketch a direct  proof by induction, as this will have singularity
conf\/inement as an immediate corollary.
The inductive hypothesis is that any four adjacent iterates
$x_k, x_{k+1},x_{k+2},x_{k+3}$ for $0\leq k\leq n+4$
are coprime elements of the unique factorization domain
$\ring =\Z [x_0^{\pm 1},x_{1}^{\pm 1},x_{2}^{\pm 1},x_{3}^{\pm 1},\be ]$,
which has units $\pm x_0^{\ell_0}, x_1^{\ell_1},x_2^{\ell_2}x_3^{\ell_3} $
for $\ell_j\in\Z$.
Working  $\bmod \,x_{n+4}$,
the congruences $x_{n+5}\equiv \be x_{n+3}^{c}x_{n+1}^{-1}$,
$x_{n+6}\equiv \be x_{n+3}^{c+1}x_{n+2}^{-1}x_{n+1}^{-1}$,
$x_{n+7}\equiv \be^{c+1} x_{n+3}^{c^2-1}x_{n+1}^{-c}$ all
hold,  so that
\[
x_{n+8}x_{n+4}\equiv x_{n+3}^{c^2+c-1}x_{n+2}^{-c}x_{n+1}^{-c-1} \,
\Big(
x_{n+3}x_{n+1}+\be\,x_{n+2}^{c}
\Big)\equiv 0,
\]
since the bracketed expression in the middle
is just $x_{n+4}x_n$ by (\ref{7betarec}). This
proves the inductive step that $x_{n+8}\in \ring$, and it is easy to see
from (\ref{7betarec}) that
this element is coprime to~$x_{n+5}$, $x_{n+6}$, $x_{n+7}$; the base of the induction is trivial.
This argument also demonstrates singularity conf\/inement: if we have
$x_{n+4}=(x_{n+3}x_{n+1}+\be\,x_{n+2}^{c})/x_n=\epsilon\to 0$
for some $n$, so that $x_{n+8}$ is potentially singular, then the preceding calculation
shows that $x_{n+8}x_{n+4}=O(\epsilon )$ and hence
$x_{n+8}=O(1)$ as $\epsilon\to 0$, so that the singularity is conf\/ined.
The Poissonicity of the map~(\ref{7m123map}) is checked by a direct calculation, and in the
coordinates $y_n=\log x_n$ the Poisson tensor
for the bracket $\{\cdot ,\cdot \}_c$ is constant and has determinant $(c-2)^2(c+1)^2$.
Thus for $c\in\N$ it is nondegenerate unless $c=2$, which gives
the previously mentioned bracket (\ref{s4brkt}) preserved by the Somos-4
recurrence.
\end{proof}

\begin{remark}
In the case of arbitrary parameters $a=b$ and $c$,  each of the recurrences (\ref{7betarec})
admits a log-canonical Poisson  bracket that generalizes $\{\cdot ,\cdot \}_c$, but
there is no such bracket for $a\neq b$.
\end{remark}

\begin{corollary}
For each $c$ the two-form
\begin{gather}
\omega    =
  \left( \frac{\dd x_0 \wedge \dd x_1}{x_0x_1}+\frac{\dd x_0 \wedge \dd x_3}{x_0x_3}+
\frac{\dd x_2 \wedge \dd x_3}{x_2x_3}\right)\nonumber \\
\phantom{\omega    =}{}
-c \left(\frac{\dd x_0 \wedge \dd x_2}{x_0x_2}+ \frac{\dd x_1 \wedge \dd x_3}{x_1x_3}\right)
+(c+1) \frac{\dd x_1 \wedge \dd x_2}{x_1x_2},\label{7sympm1is1}
\end{gather}
is preserved by the map \eqref{7m123map}, and this is symplectic for $c\neq 2$.
When $c=2$ this two-form
is degenerate, being
the pullback of the two-form $\om_1$ in \eqref{om}
under the transformation
\[
(x_0,x_1,x_2,x_3)\mapsto (u_1,u_2)= \big(x_0x_2/(x_1)^2,x_1x_3/(x_2)^2\big) .
\]
\end{corollary}

We should now like to assert that the recurrences (\ref{7betarec})
do not correspond to algebraically completely integrable systems
when $c\geq 3$,
based on the fact that in that case they have non-zero algebraic entropy.
Recall that for a rational map, the algebraic entropy is def\/ined
as $\lim\limits_{n\to\infty} (\log d_n)/n$, where $d_n$ is the degree
of the $n$th iterate of the map \cite{hv}.
Usually f\/inding  this limit requires extensive calculations of the corresponding sequence of rational
functions of the initial data, or of the iterates of the projectivized form of the map.
However, in this case we can exploit the fact that these recurrences generate Laurent
polynomials, as well as the rescaling $x_n\to \be^{ -1/(c-2)}x_n$ (which for $c \neq 2$ is equivalent
to setting $\be =1$ in  (\ref{7betarec})), to
argue that the degrees of the iterates as polynomials in the coef\/f\/icient $\be$ gives
a suitable measure of the entropy.

\begin{proposition} \label{betadegs}
For $c\in\N$, the $n$th iterate of \eqref{7betarec} is a polynomial in $\be$ of degree $d_n$,
as well as a Laurent polynomial in the initial data, where $d_n$ satisfies the
recurrence
\beq\label{7troplin}
d_{n+2}+d_{n-2}=\max \{ \, c\,d_n+1 , d_{n+1}+d_{n-1}\, \}
\eeq
for $n\geq 2$, with initial data $d_0=d_1=d_2=d_3=0$.
The  algebraic entropy of the recurrence is zero for $c=0,1,2$,
while for $c\geq 3$ it 
is given by
\beq \label{algent}
\lim_{n\to \infty} (\log d_n )/n = \frac{1}{2} \log \left(
\frac{c+\sqrt{c^2-4}}{2}
\right).
\eeq
\renewcommand{\qed}{}
\end{proposition}

\begin{remark}
The full analysis of the `tropical' (or piecewise-linear)
recurrence (\ref{7troplin}) is somewhat involved,
and is omitted here, but we can mention that the determination of the value (\ref{algent}) for
the algebraic entropy follows from the fact that when $c\geq 3$ the degrees just satisfy the
linear recurrence
\[
d_{n+2}+d_{n-2}=c\,d_n+1
\]
when $n\geq 6$, and hence they grow exponentially with $n$. For $c=0,1$ the  growth
of $d_n$ is linear in $n$, while for $c=2$ it is quadratic in $n$ (corresponding to
the quadratic growth of logarithmic heights on elliptic curves \cite{silver1}).
Very similar analysis  shows that for $c\geq 3$ the recurrences~(\ref{7m123}) fail Halburd's
Diophantine integrability criterion \cite{halburd}, which requires
that the logarithmic heights of all rational-valued iterates should grow no faster
than polynomial in $n$. For instance, with initial
data $x_0=x_1=x_2=x_3=1$ each recurrence generates polynomials in $\Z [\be ]$,
and upon evaluating these at generic values of $\be\in\Q$ it can be demonstrated that
the logarithmic heights of these numbers grow like the degrees $d_n$.
\end{remark}

Having isolated the cases $c=0,1,2$, we shall describe their integrable structure
(in descending order). The case $c=2$ is the original Somos-4
recurrence (\ref{s4}) that was treated in the previous section, so
we proceed with $c=1$.

\begin{theorem} \label{c1thm}
The map \eqref{7m123map} for $c=1$ is superintegrable, in the
sense that it has three independent conserved quantities $\J_k$, $k=1,2,3$,
which satisfy
\beq \label{jpbs}
\{ \, \J_1,\J_2\, \}_1 = 0 = \{ \, \J_1,\J_3\, \}_1,
\eeq
where
\begin{gather}
\J_1  = \Cs_0\Cs_1\Cs_2-\Cs_0^2-\Cs_1^2-\Cs_2^2+2,\qquad
\J_2  =  \Cs_0 + \Cs_1 + \Cs_2, \qquad
\J_3  =  \Cs_0\Cs_1\Cs_2.\label{jdefs}
\end{gather}
with
\begin{gather}
\Cs_0  =   \frac{x_0x_3+x_1^2+x_2^2}{x_1x_2},   \qquad
\Cs_1  =   \frac{x_0x_3^2+x_1^2x_3+x_0x_2^2+\be x_1x_2}{x_0x_2x_3},   \nonumber\\
\Cs_2  =   \frac{x_0^2x_3+x_1^2x_3+x_0x_2^2+\be x_1x_2}{x_0x_1x_3}.  \label{7cj}
\end{gather}
The iterates of the corresponding recurrence
\beq \label{7danascott}
x_{n+2}\, x_{n-2} = x_{n+1}\, x_{n-1} +  \be \, x_n,
\eeq
also satisfy the ninth-order linear recurrence
\beq \label{7ninelin}
x_{n+9}-(\J_1+1)(x_{n+6}-x_{n+3})-x_n=0,
\eeq
and the solution of the initial value problem for
\eqref{7danascott} has the form
\beq \label{7cheball}
x_{3n+j}=\rA_j\, T_n(\J_1/2)+\rB_j\,U_n(\J_1/2)+\be\,
\Cs_j/(\J_1-2),\qquad j=0,1,2,
\eeq
where $T_n$ and $U_n$
are the Chebyshev
polynomials of the first and second kind respectively, and
the coefficients $\rA_j$, $\rB_j$ are given by
\beq\label{7abj}
\rA_j=2x_j-\frac{2x_{j+3}}{\J_1}-\frac{2\be\, \Cs_j(\J_1-1)}{\J_1(\J_1-2)},
\qquad
\rB_j=-x_j+\frac{2x_{j+3}+\be\, \Cs_j}{\J_1}
\eeq
for $j=0,1,2$.
\end{theorem}

\begin{proof}
The proof of the above result is only sketched here, as further details
will be presented elsewhere \cite{ddiis}. The main observation is that
the recurrence (\ref{7danascott}) is linearizable, in the sense
that the iterates satisfy the higher-order linear recurrence (\ref{7ninelin})
for a suitable $\J_1$. (In the case $\be =1$,
the nonlinear recurrence was originally considered by Dana Scott
\cite{gale}, who found that an integer sequence was generated from the
initial data $x_0=x_1=x_2=x_3=1$; in that case the linear  recurrence (\ref{7ninelin})
is satisf\/ied with $\J_1=9$.) In general one can take (\ref{7ninelin})
as the def\/inition of $\J_1$, and use (\ref{7danascott}) to back-substitute
and rewrite it in terms of four adjacent iterates
\[
(x_n,x_{n+1},x_{n+2},x_{n+3})=(p,q,r,s),
\]
as
\beq\label{j1}
\J_1=\frac{(p^2+s^2)qr +\be (p+s)(q^2+r^2+ps)+\be^2 qr}{pqrs},
\eeq
which is found to be invariant with $n$, and def\/ines a quartic threefold in $\C^4$.
As a consequence, the linear
recurrence  (\ref{7ninelin}) holds for all $n$, and further implies the
inhomogeneous
linear equation
\begin{gather} \label{7inhgs6}
x_{n+6}-\J_1\,x_{n+3}+x_n+\be\Cs_{n} =0,
\end{gather}
where the quantity $\Cs_{n}$ varies with $n\bmod 3$. Writing everything
in terms of coordinates
\[
(x_0,x_1,x_2,x_3)\in\C^4
\]
for the map (\ref{7m123map})
with $c=1$, this gives three independent quantities $\Cs_j$
given by (\ref{7cj}) such that
\[
\{ \,\J_1,\Cs_j\, \}_1 =0, \qquad j=0,1,2.
\]
These $\Cs_j$ are not preserved by the map, but symmetric functions
of them are, which produces the formulae (\ref{jdefs}) for three
independent conserved quantities.
The solution of the initial value problem can be conveniently
expressed in the form (\ref{7cheball}), upon noting
that (by separating out $n\bmod 3$)
the homogeneous form of
(\ref{7inhgs6}) is equivalent to the second-order linear
dif\/ference equation satisf\/ied by the Chebyshev polynomials
$T_n(\J_1/2)=\cos (n\theta )$, $U_n(\J_1/2) = \sin (n\theta)/\sin\theta$
with   $\J_1 = 2\cos\theta$.
\end{proof}

\begin{remark}
The situation whereby an integrable system has more independent conserved
quantities than the number of degrees of freedom is known as
non-commutative integrability
(in~the sense of Nekhoroshev \cite{nekhor}), because not all these
quantities can be in involution with one another.
In this example,
$\J_1$ Poisson commutes with both $\J_2$ and $\J_3$, but $\{ \J_2,\J_3 \}_1 \neq 0$.
The terminology `superintegrable' is applied in the even more special
situation that the number of independent integrals is one less
than the dimension of the phase space \cite{wsuper}, as is the case~here.
\end{remark}

Upon applying the Diophantine Laurentness Lemma to the case of
initial data $(1,1,1,1)$, and choosing integer $\be$ (with $\be\neq 0$
to avoid the degenerate case of a f\/ixed point) we get
inf\/initely many solutions of certain Diophantine equations corresponding
to the conserved quantities.

\begin{corollary}\label{c1cor}
With $\J_1 =\J_1 (p,q,r,s)$ as in \eqref{j1},
\begin{gather*}
\J_2 = \frac{pr^2(q+r+s)+sq^2(p+q+r)+ps(pr+ps+qs)+\be qr(q+r)}{pqrs} ,
\\
\J_3=\frac{(ps+q^2+r^2)(pr^2+ps^2+q^2s+\be qr)(pr^2+p^2s+q^2s+\be qr)}
{p^2q^2r^2s^2} ,
\end{gather*}
there are infinitely many integer solutions $(p,q,r,s)$ of
the double pencil of Diophantine equations
given
by
\beq\label{doublepencil}
\la_1 \, \J_1 +\la_2 \, \J_2 +\la_3 \, \J_3 =
\la_1 (\be^2+6\be +2 )+\la_2 (2\be +9 )+ 3\la_3(\be +3)^2 ,
\eeq
for all $\be\in\Z\setminus \{ 0 \}$, and any
$(\la_1 :\la_2:\la_3)\in\Pro^2$.
\end{corollary}

As was remarked after Proposition \ref{betadegs},
the Laurent polynomials generated by the $c=0$ case of (\ref{7betarec}) show linear degree
growth, so it
might be anticipated that this case should also be linearizable. This indeed turns out to be so: the
main results are very similar to the case $c=1$, and are stated below without proof.

\begin{theorem} \label{c0thm}
The map \eqref{7m123map} for $c=0$ is superintegrable, in the
sense that it has three independent conserved quantities $\tiJ_j$, $j=1,2,3$,
which satisfy
\beq \label{tijpbs}
\{ \, \tiJ_1,\tiJ_2\, \}_0 = 0 = \{ \, \tiJ_1,\tiJ_3\, \}_0,
\eeq
and also $\{ \, \tiJ_1,\rQ_j\, \}_0 = 0$  for
$j=0,1,2$, where
\begin{gather}
\tiJ_1=\rQ_0\rQ_1\rQ_2-\rQ_0-\rQ_1-\rQ_2, \qquad \tiJ_2=\rQ_0\rQ_1+\rQ_1\rQ_2+\rQ_2\rQ_0-3, \nonumber\\
\tiJ_3=\rQ_0\rQ_1\rQ_2,\label{7tij123}
\end{gather}
with
\beq \label{7qs}
\rQ_0 =\frac{x_0+x_2}{x_1}, \qquad
\rQ_1 =\frac{x_1+x_3}{x_2}, \qquad
\rQ_2 =\frac{x_1}{x_0}+\frac{x_2}{x_3}+\frac{\be }{x_0x_3}.
\eeq
The iterates of the corresponding recurrence
\beq \label{7mzero}
x_{n+2}\, x_{n-2} = x_{n+1}\, x_{n-1} +  \be ,
\eeq
satisfy the sixth-order linear recurrence
\beq \label{7hgs6}
x_{n+6}-\tiJ_1\,x_{n+3}+x_n=0
\eeq
and the solution of the initial value problem for
\eqref{7mzero} can be written explicitly
in terms of Chebyshev
polynomials of the first and second kind ($T_n$ and $U_n$ respectively), as
\beq \label{7mzerocheb}
x_{3n+j}=\tilde{\rA}_j\, T_n(\tiJ_1/2)+\tilde{\rB}_j\,U_n(\tiJ_1/2),\qquad j=0,1,2,
\eeq
where
the coefficients $\tilde{\rA}_j$, $\tilde{\rB}_j$ are given by
\beq\label{7mzeroabj}
\tilde{\rA}_j=2x_j-\frac{2x_{j+3}}{\tiJ_1},
\qquad
\tilde{\rB}_j=-x_j+\frac{2x_{j+3}}{\tiJ_1}
\eeq
for $j=0,1,2$.
\end{theorem}

We can apply the Diophantine Laurentness Lemma once more, taking 
initial data $(1,1,1,1)$ and  $\be\in\Z$,
(with $\be\neq 0,-1$ to avoid periodic orbits),
to get
inf\/initely many solutions of quartic Diophantine equations corresponding
to the conserved quantities for (\ref{7mzero}).

\begin{corollary}\label{c0cor}
With the identification
$(w,x,y,z)=(x_n,x_{n+1},x_{n+2},x_{n+3})$, when $n=0$  the
relations \eqref{7tij123}  define $\tiJ_k =\tiJ_k (w,x,y,z)$ for $k=1,2,3$, as
\begin{gather*}
\tiJ_1=\frac{(w^2+z^2)xy+\beta (wx+wz+yz)}{wxyz},
\\
\tiJ_2 =\left(\bear{l} w^2z^2+xz(w^2+x^2+y^2+xz)  \\
+wy(x^2+y^2+z^2+wy)+\be (x^2+y^2+xz+wy)
\eear \right)/(
wxyz) ,
\\
\tiJ_3=\frac{(w+y)(x+z)(xz+yw+\beta )}
{wxyz} ,
\end{gather*}
and there are infinitely many integer solutions $(w,x,y,z)$ of
the double pencil of Diophantine equations 
given by
\beq\label{doublepencil0}
\la_1 \, \tiJ_1 +\la_2 \, \tiJ_2 +\la_3 \, \tiJ_3 =
\la_1 (3\be +2)+\la_2 (4\be +9 )+ 4\la_3(\be +2) ,
\eeq
for all $\be\in\Z\setminus \{ 0,-1 \}$, and any
$(\la_1 :\la_2:\la_3)\in\Pro^2$.
\end{corollary}

\begin{remark}
The initial data $(1,1,1,1)$, together with the restrictions on $\be$, are suf\/f\/icient
to ensure that the each of recurrences (\ref{7danascott}) and (\ref{7mzero}) generate
non-periodic integer sequences, and hence
inf\/initely many solutions of the corresponding Diophantine equations, given in
Corollary~\ref{c1cor} and Corollary~\ref{c0cor} respectively. However, due to the fact that the
recurrences are integrable (and even linearizable) in both cases, it is possible to choose
much more general initial data and still generate integer sequences, which produce
solutions of the same Diophantine equations but with dif\/ferent values on the right hand sides
of (\ref{doublepencil}) and (\ref{doublepencil0})
respectively.
\end{remark}

\section{Outlook}

\setcounter{equation}{0}

The Laurent property appears to be an extremely
elegant,  but somewhat special, feature of certain rational maps.
In particular, it seems to hold for integrable bilinear or discrete Hirota type
equations, such as (\ref{s4}) and (\ref{hirota}), but also for
the whole family (\ref{7betarec}),
whose members have non-zero algebraic entropy for $c\geq 3$. For the latter
family, we have noted the close connection between the Laurent property and
the notion of
singularity conf\/inement as introduced in \cite{grp}. (For
other examples of conf\/ined maps with the Laurent property
see \cite{dioph, singlaur}.)
This connection seems to persist
for rational maps that do not themselves have the Laurent property.

\begin{figure}[ht!]
\centerline{
\scalebox{1.0}{\includegraphics{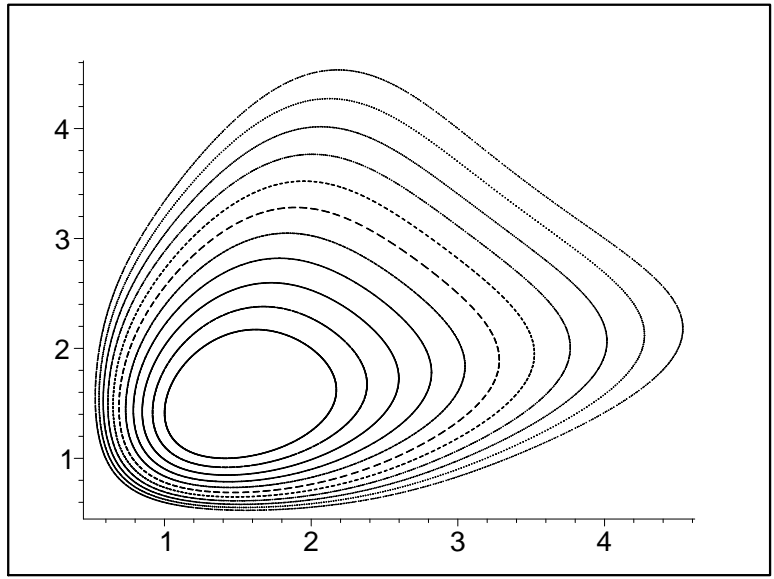}}
}
\caption{A family of orbits for the nonlinear
recurrence (\ref{nonint}).}
\label{viapaporbits}
\end{figure}

For example, consider the second-order equation\footnote{I am grateful to
Vasilis Papageorgiou for showing me this example, which I believe is due to
Claude Viallet.}
\beq \label{nonint}
u_{n+1}=\frac{(u_n)^2+1}{u_{n-1}u_n},
\eeq
which is superf\/icially very similar to (\ref{ueq}), and
preserves the same symplectic form (\ref{om}).  The real phase
portrait in $\R^2$ also looks qualitatively similar: Fig.~\ref{viapaporbits}
seems to display the same structure of invariant curves
as Fig.~\ref{s4orbits}.
Furthermore, the
equation (\ref{nonint})
satisf\/ies singularity conf\/inement, with the singularity pattern
being $\epsi$, $\epsi^{-1}$, $\epsi^{-2}$, $\epsi^{-1}$, $\epsi$ (for $\epsi\to 0$),
which suggests making the substitution
\beq
\label{taus}
u_n = \frac{\tau_{n+2}\tau_{n-2}}{\tau_{n+1}(\tau_n)^2\tau_{n-1}}.
\eeq
Thus $u_n$ given as above satisf\/ies (\ref{nonint}) whenever $\tau_n$
satisf\/ies
\beq \label{taueq}
\tau_{n+3}\, \tau_{n-3} = (\tau_{n+2}\, \tau_{n-2} )^2 + (\tau_{n+1})^2 (\tau_{n})^4(\tau_{n-1})^2 ,
\eeq
and this sixth-order recurrence has the Laurent property, as
well as satisfying the  singularity con\-f\/i\-ne\-ment test. The singularity pattern for (\ref{nonint}),
which includes poles, ``unfolds'' to yield isolated zeros, i.e. $\tau_n=\epsi$ for some $n$
with adjacent iterates being $O(1)$ as $\epsi\to 0$.

However, the logarithmic heights $h(\tau_n)$
of rational iterates grow exponentially with $n$. To see
this, it is instructive to take all six initial values for (\ref{taueq})
equal to 1, yielding the integer sequence
\[
1,1,1,1,1,1,2,5,29,1241,3642581, 80305336110269,\ldots ,
\]
which grows like $\log h(\tau_n)= \log \log |\tau_n| \sim n\log\gamma$
with $\log\gamma\approx 0.733$,
where $\gam \approx 2.081$ is the largest modulus root of the polynomial
$\gam^4-\gam^3-2\gam^2-\gam +1$. Hence the logarithmic
heights of the rational numbers $u_n$ that lie in the orbit
of $(u_0,u_1)=(1,1)$ generated by (\ref{nonint}) grow exponentially,
and Halburd's Diophantine integrability criterion is failed. Similar arguments
hold for generic orbits, and it follows that the curves appearing in Fig.~\ref{viapaporbits}
are not algebraic. To see this, recall that by the Hurwitz theorem a curve
with an inf\/inite
order automorphism group has genus zero or one, and under
iteration of such automorphisms
the logarithmic heights
of rational points grow linearly on a rational curve and quadratically
on an elliptic one \cite{silver1}.


We have concentrated on recurrences of the particular form (\ref{rec2}),
but this is not necessary for the Laurent property. Another
interesting (and algebraically non-integrable) example is the second-order
equation
\beq\label{hv}
u_{n+1}+u_{n-1}=u_n +\frac{a}{(u_n)^2}, \qquad a\neq 0,
\eeq
which in \cite{hv} was found by Hietarinta and Viallet to display singularity
conf\/inement with the pattern $\epsi$, $\epsi^{-2}$, $\epsi^{-2}$, $\epsi$,
yet it has positive algebraic entropy and its real orbits
display the characteristics of chaos.
By way of the substitution
\beq
\label{hvtaus}
u_n = \frac{\tau_{n+2}\tau_{n-1}}{(\tau_{n+1}\tau_n)^2}
\eeq
we arrive at the f\/ifth-order recurrence
\beq\label{hvfif}
\tau_{n+3}(\tau_{n})^3(\tau_{n-1})^2=(\tau_{n+2})^3(\tau_{n-1})^3-(\tau_{n+2})^2(\tau_{n+1})^3\tau_{n-2}
+a(\tau_{n+1}\tau_n)^6,
\eeq
which itself
satisf\/ies singularity conf\/inement and has the Laurent property,
i.e.\ for all $n$ the iterates satisfy
$
\tau_n\in \Z [\tau_0^{\pm 1},\ldots ,\tau_{4}^{\pm 1},a ].
$
In this case, the logarithmic heights of rational iterates $\tau_n\in\Q$
(with $a\in\Q$) generically satisfy $h(\tau_n)\sim C \ze^n$ for some $C>0$ and
$\ze = (3+\sqrt{5})/2$ is the square of the golden mean, while
$\log\ze$
turns out to be the value of the algebraic entropy for (\ref{hv}).
Note that while the calculation of the algebraic entropy is
quite involved
\cite{hv, takenawa}, it is quite straightforward to calculate
the growth of heights from (\ref{hvfif}).
Similarly to the previous example,
the only
conf\/ined singularities
that appear in (\ref{hvfif}) are isolated zeros.

These examples illustrate the following general phenomenon: whenever
we have a rational map with conf\/ined singularities, including poles,
it should always be possible
to ``unfold'' these into conf\/ined zeros, by embedding the map in
higher dimensions via a change of variables, and the new map
thus obtained should have the
Laurent property. This is analogous to the way in which continuous
integrable systems with the Painlev\'e property, that have meromorphic solutions,
admit a Hirota bilinear form
(or multilinear form) in terms of tau-functions that are holomorphic.
Although this phenomenon (the existence of a tau-function)
is very well known for discrete integrable systems \cite{bilins}, so far it does
not seem to have been to have been exploited in the case
of non-integrable maps.

\pdfbookmark[1]{References}{ref}
\LastPageEnding


\begin{thebibliography}{99}

\footnotesize\itemsep=0pt

\bibitem {beh} Braden H.W., Enolskii V.Z.,  Hone A.N.W.,
Bilinear recurrences and addition formulae for hyperelliptic sigma
functions,
{\it J. Nonlinear Math. Phys.} {\bf 12} (2005),
suppl.~2, 46--62, \href{http://arxiv.org/abs/math.NT/0501162}{math.NT/0501162}.

\bibitem{bressoud}
Bressoud D.M.,
Proofs and conf\/irmations: the story of the alternating sign matrix
conjecture, Cambridge University Press, Cambridge, 1999.

\bibitem {rag} Bruschi M., Ragnisco O., Santini P.M.,
Tu G.-Z.,
Integrable symplectic maps,
{\it Phys.~D} {\bf 49} (1991), 273--294.

\bibitem{heron}
Buchholz R.H., Rathbun R.L.,
An inf\/inite set of heron triangles
with two rational medians,
{\it Amer. Math. Monthly} {\bf 104} (1997), 107--115.

\bibitem{bel}Buchstaber V.M., Enolskii V.Z., Leykin D.V.,
Hyperelliptic Kleinian functions and applications,
in Solitons, Geometry and Topology: On the Crossroad,
Editors V.M.~Buchstaber and S.P.~Novikov,
{\it AMS Translations Series 2}, Vol.~179, AMS,
1997, 1--34, \href{http://arxiv.org/abs/solv-int/9603005}{solv-int/9603005}.

\bibitem{cantor} Cantor D.,
On the analogue of
the division polynomials for hyperelliptic curves,
{\it J. Reine Angew. Math.} {\bf 447} (1994),
91--145.

\bibitem{cube}Carroll G., Speyer D.,
The cube recurrence,
{\it Electron. J. Combin.} {\bf 11} (2004),
$\#\,$R73, \href{http://arxiv.org/abs/math.CO/0403417}{math.CO/0403417}.

\bibitem{chm}Common A., Hone A.N.W., Musette M., A new discrete H\'enon--Heiles system,
{\it J. Nonlinear Math. Phys.} {\bf 10} (2003), suppl.~2, 27--40.

\bibitem{chaos} Elaydi S.,
Discrete chaos, Chapman and Hall/CRC,
Boca Raton, 2000.

\bibitem {eew} Einsiedler M., Everest G., Ward T.,
Primes in elliptic divisibility sequences,
{\it LMS J. Comput. Math.} {\bf 4} (2001),
1--13.


\bibitem {ems} Everest G., Miller V., Stephens N.,
Primes generated by elliptic curves,
{\it Proc. Amer. Math. Soc.} {\bf 132} (2003),
955--963.

\bibitem {recs} Everest G., van der Poorten A.,
Shparlinski I., Ward T., Recurrence sequences,
{\it AMS  Mathematical Surveys and Monographs}, Vol.~104,
 Amer. Math. Soc., Providence, RI,   2003.

\bibitem{fed2}
Fedorov Y.,
B\"{a}cklund transformations on
coadjoint orbits of the loop algebra $gl(r)$,
{\it J. Nonlinear Math. Phys.} {\bf 9} (2002),
suppl.~1,  29--46.

\bibitem {fz} Fomin S., Zelevinsky A.,
The Laurent phenomenon,
{\it Adv. Appl. Math.} {\bf 28} (2002),  119--144, \href{http://arxiv.org/abs/math.CO/0104241}{math.CO/0104241}.

\bibitem {fz4} Fomin S., Zelevinsky A.,
Cluster algebras IV: coef\/f\/icients,
{\it Compos.  Math.},  to appear,
\href{http://arxiv.org/abs/math.RA/0602259}{math.RA/0602259}.



\bibitem{poissonmaps} Fordy A.P., Shabat A.B., Veselov A.P.,
Factorization and Poisson correspondences,
{\it Theor. Math. Phys.} {\bf 105} (1995), 1369--1386.

\bibitem {gale} Gale D.,
The strange and surprising saga of the
Somos sequences,
{\it Math. Intelligencer} {\bf 13} (1991), no.~1, 40--42.\\
Gale D., Somos sequence update,
{\it Math. Intelligencer} {\bf 13} (1991), no.~4,
49--50 (Reprinted in Tracking the Automatic Ant.,
Springer, 1998).

\bibitem{gekhtman}
Gekhtman M., Shapiro M., Vainshtein A.,
Cluster algebras and Poisson geometry,
{\it Moscow Math. J.} {\bf 3} (2003),
899--934, \href{http://arxiv.org/abs/math.QA/0208033}{math.QA/0208033}.


\bibitem{grp} Grammaticos B., Ramani A., Papageorgiou V.,
Do integrable mappings have the Painlev\'e property?
{\it Phys. Rev. Lett.} {\bf 67} (1991),
1825--1828.

\bibitem{halburd} Halburd R.G., Diophantine integrability,
{\it J. Phys. A: Math. Gen.} {\bf 38} (2005),  L263--L269,
\href{http://arxiv.org/abs/nlin.SI/0504027}{nlin.SI/0504027}.

\bibitem{hv} Hietarinta J., Viallet C.,
Singularity conf\/inement and chaos in discrete systems,
{\it Phys. Rev. Lett.}  {\bf 81} (1998), 325--328, \href{http://arxiv.org/abs/solv-int/9711014}{solv-int/9711014}.

\bibitem{nahh} Hone A.N.W.,
Non-autonomous H\'enon--Heiles systems,
{\it Phys.~D} {\bf 118} (1998), 1--16, \href{http://arxiv.org/abs/solv-int/9703005}{solv-int/9703005}.

\bibitem{hkr1} Hone A.N.W.,
Kuznetsov V.B., Ragnisco O.,
B\"acklund transformations for the H\'enon--Heiles and Garnier systems,
{\it CRM Proceedings and Lecture Notes}, Vol.~25, Amer. Math. Soc.,  2000,
231--235.

\bibitem{hkr2} Hone A.N.W.,
Kuznetsov V.B., Ragnisco O.,
B\"acklund transformations for many-body systems related to KdV,
{\it J. Phys. A: Math. Gen.} {\bf 32} (1999), L299--L306, \href{http://arxiv.org/abs/solv-int/9904003}{solv-int/9904003}.


\bibitem{hkr3} Hone A.N.W.,
Kuznetsov V.B., Ragnisco O.,
B\a"{a}cklund transformations for the $sl(2)$ Gaudin magnet,
{\it J.~Phys.~A: Math. Gen.} {\bf 34} (2001), 2477--2490, \href{http://arxiv.org/abs/nlin.SI/0007041}{nlin.SI/0007041}.

\bibitem{pinney} Hone A.N.W.,
Exact discretization of the Ermakov--Pinney equation,
{\it Phys. Lett. A} {\bf 263} (1999), 347--354.


\bibitem {honeblms} Hone A.N.W.,
Elliptic curves and quadratic recurrence sequences,
{\it Bull. Lond. Math. Soc.} {\bf 37} (2005), 161--171,
Corrigendum, {\it Bull. Lond. Math. Soc.},  {\bf 38} (2006), 741--742.

\bibitem{hones5} Hone A.N.W.,
Sigma function solution of the initial value problem
for Somos 5 sequences,
{\it Trans. Amer. Math. Soc.}, to appear,
\href{http://arxiv.org/abs/math-ph/0312029}{math.NT/0501554}.

\bibitem{dioph} Hone A.N.W.,
Diophantine non-integrability of
a third-order recurrence with the Laurent property,
{\it J.~Phys.~A: Math. Gen.} {\bf 39} (2006), L171--L177, \href{http://arxiv.org/abs/math.NT/0601324}{math.NT/0601324}.

\bibitem{singlaur} Hone A.N.W.,
Singularity conf\/inement for maps with the Laurent property,
{\it Phys. Lett. A} {\bf 361} (2007), 341--345, \href{http://arxiv.org/abs/nlin.SI/0602007}{nlin.SI/0602007}.

\bibitem{ddiis} Hone A.N.W.,
Discrete dynamics, integrability and integer sequences,
Imperial College Press, in preparation.

\bibitem{kuskly}Kuznetsov V.B., Sklyanin E.K.,
On B\a"{a}cklund transformations for many-body systems,
{\it J. Phys. A: Math. Gen.} {\bf 31} (1998), 2241--2251, \href{http://arxiv.org/abs/solv-int/9711010}{solv-int/9711010}.

\bibitem{kuskly2}Kuznetsov V.B., Salerno M., Sklyanin E.K.,
Quantum B\a"{a}cklund transformation for the integrable DST model,
{\it J. Phys. A: Math. Gen.} {\bf 33} (2000), 171--189, \href{http://arxiv.org/abs/solv-int/9908002}{solv-int/9908002}.

\bibitem{vadimpol}
Kuznetsov V.B., Vanhaecke P.,
B\a"{a}cklund transformations for f\/inite-dimensional integrable systems:
a geometric approach,
{\it J. Geom. Phys.} {\bf 44} (2002),
1--40, \href{http://arxiv.org/abs/nlin.SI/0004003}{nlin.SI/0004003}.

\bibitem{qopsjack}Kuznetsov V.B., Mangazeev V.V., Sklyanin E.K.,
$Q$-operator and factorised separation chain for Jack polynomials,
{\it Indag. Math.} {\bf 14} (2003), 451--482, \href{http://arxiv.org/abs/math.CA/0306242}{math.CA/0306242}.

\bibitem{mrr}
Mills W.H.,  Robbins D.P., Rumsey H.,
Alternating-sign matrices and descending plane partitions,
{\it J. Combin. Theory Ser. A} {\bf 34} (1983), 340--359.

\bibitem{nekhor}Nekhoroshev N.N.,
On action-angle variables and their generalizations,
{\it Tr. Moscow Math. Soc.} {\bf 26} (1972), 181--198 (in Russian).

\bibitem{pg}Pasquier V., Gaudin M.,
The periodic Toda chain and a matrix generalization of the Bessel function
recursion relations,
{\it J. Phys. A: Math. Gen.} {\bf 25} (1992), 5243--5252.

\bibitem{kupera}
Kuznetsov V.B., Petrera M., Ragnisco O.,
Separation of variables and
B\a"{a}cklund transformations for the symmetric Lagrange top,
{\it J. Phys. A: Math. Gen.} {\bf 37} (2004), 8495--8512, \href{http://arxiv.org/abs/nlin.SI/0403028}{nlin.SI/0403028}.

\bibitem{matsupsi}Matsutani S.,
Recursion relation of hyperelliptic PSI-functions
of genus two,
{\it Int. Transforms Spec. Func.} {\bf 14} (2003), 517--527, \href{http://arxiv.org/abs/math-ph/0105031}{math-ph/0105031}.


\bibitem {vdp} van der Poorten A.J.,
Elliptic curves and continued fractions,
{\it J. Integer Sequences} {\bf 8} (2005), Article 05.2.5, 19 pages,
\href{http://arxiv.org/abs/math.NT/0403225}{math.NT/0403225}.


\bibitem {swartvdp} van der Poorten A.J., Swart C.S.,
Recurrence relations for elliptic sequences: every Somos 4
is a Somos $k$,
{\it Bull. Lond. Math. Soc.} {\bf 38} (2006), 546--554, \href{http://arxiv.org/abs/math.NT/0412293}{math.NT/0412293}.

\bibitem{propp}Propp J.,
The many faces of alternating-sign matrices,
{\it Disc. Math. Theoret.
Comp. Sci. Proc.} \textbf{AA (DM-CCG)} (2001),  43--58, \href{http://arxiv.org/abs/math.CO/0208125}{math.CO/0208125}.

\bibitem{proppweb}Propp J., The ``bilinear'' forum,
\url{http://www.math.wisc.edu~propp/}.

\bibitem {qrt} Quispel G.R.W., Roberts J.A.G., Thompson C.J.,
Integrable mappings and soliton equations II,
{\it Phys.~D} {\bf 34} (1989),
183--192.

\bibitem{bilins}
Ramani A., Grammaticos B., Satsuma J.,
Bilinear discrete Painlev\a'{e} equations,
{\it J. Phys. A: Math. Gen.} {\bf 28} (1995), 4655--4665.

\bibitem {rob} Robinson R.,
Periodicity of Somos sequences,
{\it Proc. Amer. Math. Soc.} {\bf 116} (1992), 613--619.

\bibitem{silver1} Silverman J.H.,
The arithmetic of elliptic curves,
Springer,  1986.

\bibitem {silverman} Silverman J.H.,
$p$-adic properties of division polynomials
and elliptic divisibility sequences,
{\it Math. Annal.} {\bf 332} (2005),
443--471, Addendum, 473--474, \href{http://arxiv.org/abs/math.NT/0404412}{math.NT/0404412}.

\bibitem{sloane} Sloane N.J.A.,
On-line encyclopedia of integer sequences,
\mbox{\url{http://www.research.att.com/~njas/sequences}},
sequence A006720.

\bibitem{speyer} Speyer D.,
Perfect matchings and the octahedron recurrence, 2004,
\href{http://arxiv.org/abs/math.CO/0402452}{math.CO/0402452}.


\bibitem{surisbook}
Suris Y.B.,
The problem of integrable discretization: Hamiltonian approach,
{\it Progress in Mathematics}, Vol.~219, Birkh\"auser, Basel, 2003.

\bibitem{swart} Swart C.S.,
Elliptic curves and related sequences,
PhD thesis, Royal Holloway, University of London, 2003.

\bibitem{swahon}Swart C.S., Hone A.N.W.,
Integrality and the Laurent phenomenon for Somos 4
sequences,
\href{http://arxiv.org/abs/math.NT/0508094}{math.NT/0508094}.

\bibitem{takenawa}Takenawa T.,
A geometric approach to singularity conf\/inement and
algebraic entropy,
{\it J. Phys. A: Math. Gen.} {\bf 34} (2001),
L95--L102, \href{http://arxiv.org/abs/nlin.SI/0011037}{nlin.SI/0011037}.

\bibitem{vanhaecke}Vanhaecke P.,
Integrable systems in the realm of algebraic geometry,
2nd ed., Springer, 2005.

\bibitem{ves1}Veselov A.P.,
Integrable maps,
{\it Russ. Math. Surveys} {\bf 46} (1991),
1--51.

\bibitem{ves2}Veselov A.P.,
What is an integrable mapping?
in What is Integrability? Editor  V.E.~Zakharov,
Springer-Verlag, 1991, 251--272.

\bibitem{ves3}Veselov A.P., Growth and integrability in the dynamics of mappings,
{\it Comm. Math. Phys.} {\bf 145} (1992),
181--193.

\bibitem{dressing}
Veselov A.P.,
Shabat A.B.,
A dressing chain and the
spectral theory of the Schr\"odinger operator,
{\it Funct. Anal.
Appl.} {\bf 27} (1993),
1--21.



\bibitem {ward1} Ward M.,
Memoir on elliptic divisibility
sequences,
{\it Amer. J. Math.} {\bf 70} (1948),
31--74.

\bibitem {ward2} Ward M.,
The law of repetition of primes in an elliptic divisibility
sequence,
{\it Duke Math. J.} {\bf 15} (1948),
941--946.

\bibitem{weiss}
Weiss J.,
Periodic f\/ixed points of B\a"{a}cklund
transformations and the Korteweg--de Vries equation,
{\it J. Math. Phys.} {\bf 27} (1986),
2647--2656.

\bibitem{wsuper}Wojciechowski S.,
Superintegrability of the Calogero--Moser system,
{\it Phys. Lett. A} {\bf 95} (1983),
279--281.


\bibitem{yka} Yuzbashyan E.A., Kuznetsov V.B.,  Altshuler B.L.,
Integrable dynamics of coupled Fermi--Bose condensates,
{\it Phys. Rev. B} {\bf 72} (2005),
144524, 9 pages, \href{http://arxiv.org/abs/cond-mat/0506782}{cond-mat/0506782}.

\bibitem{yake} Yuzbashyan E.A., Altshuler B.L.,
Kuznetsov V.B.,
Enolskii V.Z.,
Nonequilibrium Cooper pairing in the nonadiabatic regime,
{\it Phys. Rev. B} {\bf 72} (2005),
220503(R), 4 pages, \href{http://arxiv.org/abs/cond-mat/0505493}{cond-mat/0505493}.


\bibitem{zagier}Zagier D., Problems posed at the St. Andrews
Colloquium, 1996, Solutions, 5th day, available at
\mbox{\url{http://www-groups.dcs.st-and.ac.uk/~john/Zagier/Problems.html}}.

\bibitem{zabrodin}
Zabrodin A., A survey of Hirota's dif\/ference equations,
{\it Teor. Mat. Fiz.} {\bf 113} (1997), 179--230 (in Russian),
\mbox{\href{http://arxiv.org/abs/solv-int/9704001}{solv-int/9704001}}.

\end{thebibliography}
\end{document}